\documentclass[runningheads]{llncs}
\usepackage[T1]{fontenc}
\usepackage{graphicx}
\usepackage{microtype}
\usepackage{graphicx}
\usepackage{wrapfig}
\usepackage{subfigure}
\usepackage{booktabs} 
\usepackage{multirow}
\usepackage{stfloats}

\usepackage[utf8]{inputenc} 
\usepackage[T1]{fontenc}    
\usepackage{hyperref}       
\usepackage{url}            
\usepackage{booktabs}       
\usepackage{amsfonts}       
\usepackage{nicefrac}       
\usepackage{microtype}      
\usepackage{xcolor}         

\usepackage{amsthm}
\usepackage{amsmath}
\usepackage{amssymb}
\usepackage{mathtools}
\usepackage{amsthm}
\usepackage{color}
\usepackage{csquotes}
\usepackage{multirow}
\RequirePackage{algorithm}
\RequirePackage{algorithmic}

\usepackage{comment}

\theoremstyle{plain}

\theoremstyle{definition}

\begin{document}
\title{Mixed-Integer Linear Optimization via Learning-Based Two-Layer Large Neighborhood Search}
\titlerunning{MILP via learning-based TLNS}
%

\author{%
  Wenbo Liu\inst{1}
  \and
  Akang Wang\inst{2,3}
  \and
  Wenguo Yang\inst{1}
  \and
  Qingjiang Shi\inst{2,4}
}
\authorrunning{W. Liu et al.}
\institute{University of Chinese Academy of Sciences, Beijing, China \and Shenzhen Research Institute of Big Data, China \and The Chinese University of Hong Kong, Shenzhen, China \and Tongji University, Shanghai, China}
\maketitle
\begin{abstract}
Mixed-integer linear programs (MILPs) are extensively used to model practical problems such as planning and scheduling. A prominent method for solving MILPs is large neighborhood search (LNS), which iteratively seeks improved solutions within specific neighborhoods. Recent advancements have integrated machine learning techniques into LNS to guide the construction of these neighborhoods effectively. However, for large-scale MILPs, the search step in LNS becomes a computational bottleneck, relying on off-the-shelf solvers to optimize auxiliary MILPs of substantial size. To address this challenge, we introduce a two-layer LNS (TLNS) approach that employs LNS to solve both the original MILP and its auxiliary MILPs, necessitating the optimization of only small-sized MILPs using off-the-shelf solvers. Additionally, we incorporate a lightweight graph transformer model to inform neighborhood design. We conduct extensive computational experiments using public benchmarks. The results indicate that our learning-based TLNS approach achieves remarkable performance gains--up to $66\%$ and $96\%$ over LNS and state-of-the-art MILP solvers, respectively.
\end{abstract}
\section{Introduction} 
\label{sec:Intro}

\textit{Mixed-integer linear programs} (MILPs) have become a cornerstone in various industrial applications, including network design~\cite{luathep2011global}, production planning~\cite{pochet2006production}, and route optimization~\cite{toth2002vehicle}. 
The resolution of MILPs typically poses $\mathcal{NP}$-hard challenges, with general-purpose MILP solvers resorting to the branch-and-bound method for systematic enumeration of candidate solutions. 
However, tackling large-scale MILPs with a branch-and-bound algorithm proves computationally demanding due to its exhaustive search nature. 
In practical scenarios, primal heuristic methods are commonly employed to efficiently identify high-quality feasible solutions. Although these heuristics do not guarantee optimality, they consistently deliver outstanding solutions for significantly larger MILPs.

One of the most prominent heuristic methods for addressing MILPs is \textit{large neighborhood search} (LNS)~\cite{shaw1998using}.
LNS often refines an incumbent solution by iteratively constructing a neighborhood of interest and searching within such a region via optimization of auxiliary MILPs.
Lots of efforts~\cite{berthold2014rens,danna2005exploring,localbranching,huang2023local,rothberg2007evolutionary} have been devoted to building neighborhoods by use of heuristics, delegating the search step to off-the-shelf MILP solvers.

Recently, \textit{machine learning} (ML) techniques have been extensively utilized to expedite the optimization of MILPs~\cite{ml4co,nair2020solving,zhang2023survey}.
These endeavors have been prompted by the recognition that MILPs arising from similar applications often exhibit recurrent patterns, which can be effectively captured through ML techniques.
LNS also benefits from ML techniques~\cite{huang2023searching,song2020general,sonnerat2021learningnew,wu2021learning}, where neighborhood construction is informed by utilizing \textit{graph neural networks}~(GNNs) on the graph representation of MILPs.
Despite promising advancements, these learning-based LNS methods still have ample room for improvement.
On one hand, while neighborhood construction shows potential for enhancement, it often overlooks improvements in the search step.
On the other hand, the number of GNN layers (typically $2$) is limited due to issues with over-smoothing \cite{li2018deeper}, rendering these applications incapable of effective message passing between distant variable nodes.

\begin{figure}[t]
\centering
\centerline{\includegraphics[width=\textwidth]{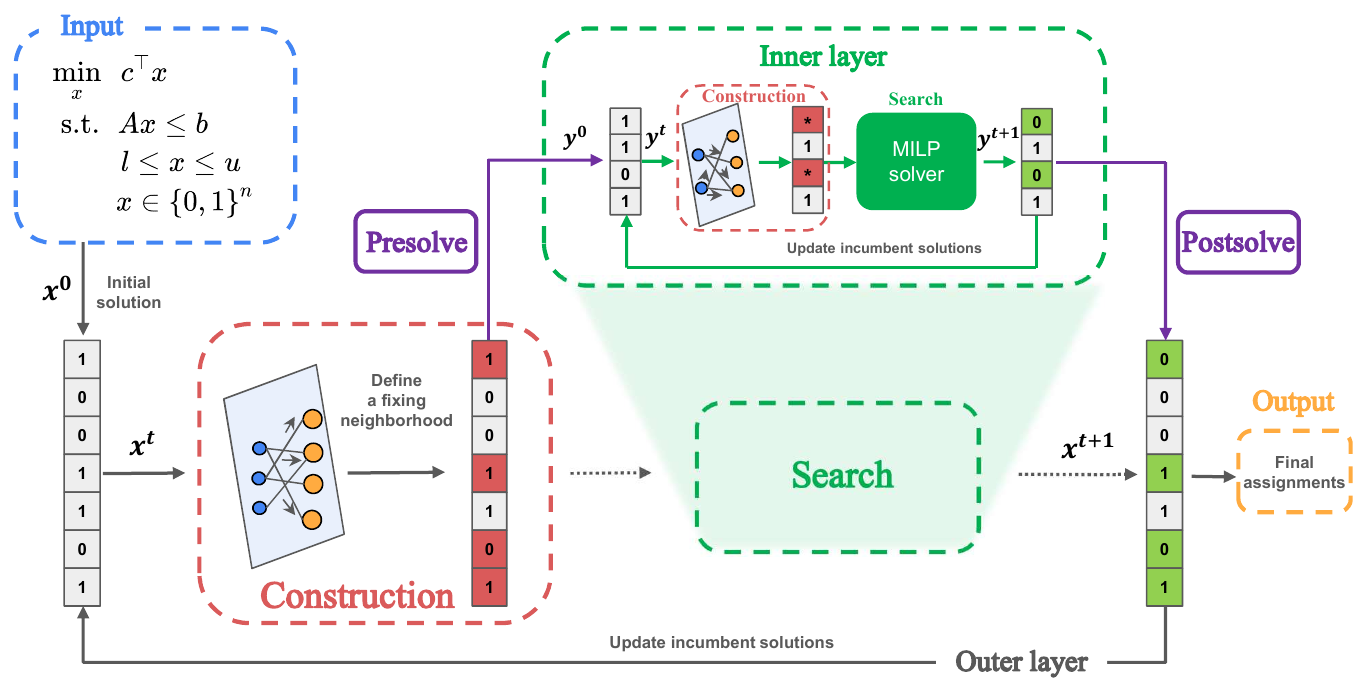}}
\caption{An overview of our proposed learning-enhanced TLNS framework. 
The red part represents the learning-enhanced neighborhood construction stage while the green part denotes the neighborhood search stage. 
The purple part indicates that the presolve operator transforms the original MILP to a reduced MILP while the postsolve operator reverses the transformation.}
\label{alg_fig}
\end{figure}
Previous efforts, in both classic and learning-enhanced settings, focused on designing effective neighborhoods, leaving the search procedure to off-the-shelf MILP solvers. 
However, relying on exact solvers to optimize auxiliary MILPs during the search stage could still be computationally expensive for large-sized problems.
Recognizing that searching for high-quality solutions still involves optimizing MILPs, we adopt a learning-based LNS approach once more. In our method, we propose applying learning-based LNS to tackle both the original MILP and its auxiliary MILPs, which necessitates optimizing only small-sized sub-MILPs using off-the-shelf solvers. Additionally, we employ a lightweight \textit{graph transformer} model to expand the receptive field of GNNs.
We call this method \enquote{learning-based Two-Layer LNS} (TLNS).
The overall algorithm is outlined in~\ref{alg_fig}.
The distinct contributions of our work can be summarized as follows.
\begin{itemize}
    \item 
    We introduce a novel TLNS algorithm designed to identify high-quality solutions for MILPs. 
    Unlike traditional approaches, TLNS applies LNS to optimize not only the original MILP (outer layer) but also its auxiliary MILPs (inner layer), employing a divide-and-conquer strategy. 
    To the best of our knowledge, this represents the first attempt to extend LNS to a multi-layer version specifically tailored for addressing MILPs.
    
    \item We employ a lightweight graph transformer model trained with \textit{contrastive loss} to effectively guide the construction of neighborhoods in LNS, furthermore boosting the algorithmic performance.
    
    \item We conduct extensive computational experiments on public benchmarks and the results show that TLNS achieves up to $66\%$ and $96\%$ improvements over LNS and the state-of-the-art MILP solvers, respectively.
\end{itemize}
The remainder of this paper is organized as follows. 
Section~\ref{sec:preliminary} introduces preliminaries along with basic techniques that will be used in our method. 
Section~\ref{sec:TLNS} presents our proposed TLNS method and how it is enhanced by ML techniques.
In Section~\ref{sec:experiments}, we conduct experiments to evaluate the performance of our method and finally we draw conclusions and discuss future work in Section~\ref{sec:conclusion}.

\subsection{Related Works}
\subsubsection{Traditional LNS}~\cite{rothberg2007evolutionary} proposed to use a mutation neighborhood by fixing a random subset of integer variables at the incumbent.
\cite{localbranching} introduced the \textit{local branching} (LB) method, defining a neighborhood as a Hamming ball around the incumbent solution, while \cite{huang2023local} proposed to utilize solutions to continuous relaxations in LB for building neighborhoods.
RINS~\cite{danna2005exploring}  constructs a promising neighborhood using information contained in the continuous relaxation of an
MILP as well as the incumbent, while RENS~\cite{berthold2014rens} relies purely on the continuous relaxation.
Though many aforementioned LNS heuristics have been deployed within MILP solvers, their computational effort makes it impractical to apply all of them frequently.
One exception is the work of~\cite{ALNS} in which the author considered eight popular LNS heuristics and proposed to adaptively select one of them for execution via online learning.

\subsubsection{Learning-based LNS} The first attempt to enhance LNS with ML is~\cite{song2020general}, in which the authors proposed to imitate the best neighborhood out of a few randomly sampled ones.
Building upon this, \cite{sonnerat2021learningnew} improved the imitation learning approach by employing LB as an expert.
Furthermore, \cite{huang2023searching} collected both positive and negative solution samples from LB and then utilized contrastive loss to learn and construct neighborhoods. 
Alternatively, \textit{reinforcement learning} (RL) was applied by~\cite{wu2021learning} to explore a promising policy for constructing neighborhoods while~\cite{liu2022learning} focused on enhancing LB and utilized RL to inform the radius of the Hamming ball.

\section{Preliminaries}
\label{sec:preliminary}
\subsubsection{Mixed-Integer Linear Programs}
An MILP is formulated as :
\begin{equation}
\setlength\abovedisplayskip{1pt}
\setlength\belowdisplayskip{0.5pt}
    \begin{aligned}
          &  \underset{x\in S}{\min} && c^\top x 
    \end{aligned} \label{eq:MILP}
\end{equation}
where $x$ denotes the decision variable and $S \coloneqq \left\{x\in \mathbb Z^q\times \mathbb R^{n-q}: Ax \leq b,l \leq x \leq u \right\}$ represents the feasible region for $x$. 
$l, u, c \in \mathbb{R}^n, b \in \mathbb{R}^m$ and $A \in \mathbb{R}^{m\times n}$ are given parameters.
For the sake of simplicity, we assume that all variables are binary, i.e., $x_i \in \left\{0, 1\right\} \; \forall i = 1, 2, ..., n$.
Let $M := (A, b, c, l, u, q)$ denote an MILP instance for convenience.

\subsubsection{Large Neighborhood Search}
In this work, we focus on \textit{fixing neighborhood LNS heuristics}~\cite{ALNS} in which neighborhoods are defined by fixing part of decision variables.

\begin{definition}[Fixing neighborhood]
\label{def:fixing neighborhood}
Consider an MILP with $n$ variables, let $\mathcal F \subseteq \left\{1,...,n\right\}$ denote an index set and $\bar{x}$ denote a reference point.
Then a \textit{fixing neighborhood} of $\bar{x}$ is defined by fixing variables in $\mathcal{F}$ to their values in $\bar{x}$:
$\mathcal{B}\left(\mathcal{F}, \bar{x} \right)\coloneqq \left\{x\in \mathbb R^n : x_i = \bar{x}_i, \forall i\in \mathcal F\right\}.$
\end{definition}
The number of unfixed variables (aka \textit{neighborhood size}) is equal to $n - |\mathcal F|$. 
Using such a fixing neighborhood, we can define an auxiliary problem.
\begin{definition}[Auxiliary problem]
Given an MILP $M$ of form~(\ref{eq:MILP}) and a fixing neighborhood $\mathcal{B}\left(\mathcal{F}, \bar{x}\right)$, then an auxiliary problem $\mathcal A(M, \bar x, \mathcal F)$ is defined as the following MILP:
\begin{equation}
\setlength\abovedisplayskip{1pt}
\setlength\belowdisplayskip{0pt}
    \begin{aligned}
          &  \underset{x \in S \cap \mathcal{B}\left(\mathcal{F}, \bar{x}\right)}{\min} && c^\top x.
    \end{aligned} \label{eq:model_auxi}
\end{equation}
\end{definition}
Problem~(\ref{eq:model_auxi}) guarantees feasibility of returned solutions whenever $\bar{x} \in S$.
In LNS, an incumbent $\bar{x}$ is consistently considered as the reference point and will be iteratively refined by constructing a fixing neighborhood and invoking an off-the-shelf solver to optimize the auxiliary problem~(\ref{eq:model_auxi}).

\subsubsection{Bipartite Graph Representation}
\begin{figure}[h]
\label{fig:BG}
\centering
\centerline{\includegraphics[width=0.7\textwidth]{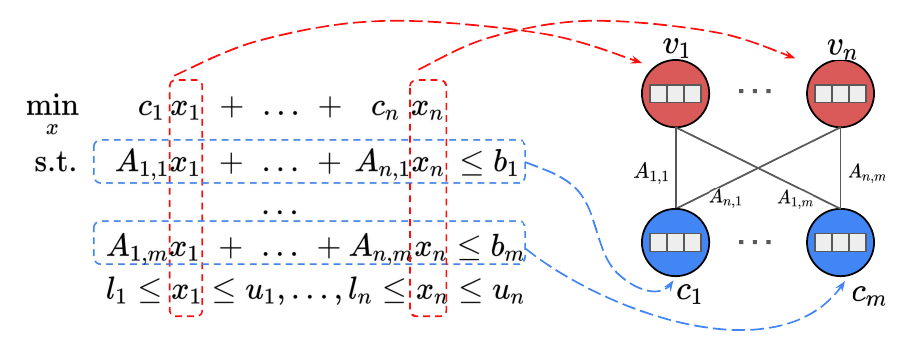}}
\caption{The bipartite graph representation of an MILP, where node $v_i$ and $v_{n+j}$ indicates the $i$-th variable and the $j$-th constraint, respectively.}
\label{figure:bipartite}

\vskip -0.2in
\end{figure}

Given an MILP of form~(\ref{eq:MILP}), \cite{gasse2019exact} proposed a variable-constraint bipartite graph representation, as shown in~\ref{figure:bipartite}.
Specifically, let $G \coloneqq \left(\mathcal{V}, \mathcal{E}\right)$ denote a bipartite graph, where $\mathcal{V} \coloneqq \left\{v_1, ..., v_n, v_{n+1}, ..., v_{n+m}\right\}$ denotes the set of $n$ variable nodes and $m$ constraint nodes, and $\mathcal{E}$ represents the set of edges that only connect between nodes of different types. 
Variable node $v_i$ and constraint node $v_{n+j}$ are connected if $A_{ji}$ is non-zero.
The information of an MILP including $c, l, u$, $b$ and $A$ will be properly incorporated into $G$ as graph attributes. 

\subsubsection{Graph Neural Networks}
For a graph $G =(\mathcal V,\mathcal E)$, let $\mathcal N(v)$ denote the set of neighbors of $v$.
The $k$-th message passing layer updates embeddings for each node $v$ using the following formula:
\begin{equation*}
\setlength\abovedisplayskip{2.5pt}
\setlength\belowdisplayskip{0.5pt}
    h_v^{(k)}=f_2^{(k)} \left( \left\{ h_v^{(k-1)},f_1^{(k)} \left ( \left\{h_u^{(k-1)}: u\in \mathcal N(v)  \right\}  \right) \right\} \right),
\end{equation*}
where $h_v^{(k)}\in \mathbb R^d$ denotes the hidden feature vector of node $v$ in the $k$-th layer with $h_v^{(0)}$ being the initial embedding.
Function $f_1^{(k)}(\cdot)$ is the \textit{AGGREGATE} operator that gathers information from neighbors while function $f_2^{(k)}(\cdot)$ is the \textit{COMBINE} operator that updates the aggregated representation.
These two operators can take various choices, resulting in different architectures such as \textit{graph convolutional networks}~\cite{gcn} and \textit{graph attention networks} (GATs)~\cite{gat}.

\subsubsection{Simplified Graph Transformers} 
Transformers with global attention~\cite{vaswani2017attention} can be considered a generalization of message passing to a fully connected graph.
Typically,~\cite{wu2024simplifying} proposed a simplified graph transformer that incorporates GNNs with a \textit{linear attention function} defined as follows:
\begin{equation}
\setlength\abovedisplayskip{1pt}
\setlength\belowdisplayskip{0pt}
    \mathbf Q = f_Q(\mathbf H^{(0)}),\quad \tilde {\mathbf Q} = \frac{\mathbf Q}{\|\mathbf Q\|_\mathcal F}, \quad \mathbf K = f_K(\mathbf H^{(0)}), \quad \tilde {\mathbf K} = \frac{\mathbf K}{\|\mathbf K\|_\mathcal F}, \quad \mathbf V = f_V(\mathbf H^{(0)}), 
\end{equation}
\begin{equation}
\setlength\abovedisplayskip{1pt}
\setlength\belowdisplayskip{0pt}
    \mathbf D = \text{diag}\left(\mathbf{1}+\frac1N\tilde{ \mathbf Q}(\tilde{\mathbf K}^\top \mathbf{1} )\right), \quad \mathbf H = \beta \mathbf D^{-1}\left[ \mathbf V+\frac1N\tilde {\mathbf Q}(\tilde {\mathbf K}^\top \mathbf V)\right]+(1-\beta)\mathbf H^{(0)}\label{eq:sgt}
\end{equation}
where $\mathbf H^{(0)}\in \mathbb R^{|\mathcal V|\times d}$ represents the initial node embeddings, and $f_Q, f_K, f_V$ denote shallow neural layers.
$\mathbf H$ is the output of the attention module with $\beta$ serving as a hyper-parameter for residual link.
Given the simplified attention module,~\cite{wu2024simplifying} utilized GNNs to incorporate structural information by adding the outputs of the two modules: $\mathbf H_O = (1-\alpha)\mathbf H +\alpha \text{GNN}(\mathbf H^{(0)})$.


\section{The Learning-Based Two-Layer Large Neighborhood Search}
\label{sec:TLNS}
In this section, we will provide a detailed description of the TLNS method and discuss how machine learning techniques contribute to fixing neighborhoods in TLNS.

\subsection{A TLNS Framework}
\begin{minipage}{0.486\textwidth}
\begin{algorithm}[H]
\caption{Large Neighborhood Search (LNS)}
\begin{algorithmic}[1]
\label{alg:LNS}
\STATE {\bfseries Input:} an MILP $M$,
   initial solution $\bar{x}$, fixing heuristic $\mathcal D$, count limit $C$, neighborhood size $k$, and adaptive rate $\eta$.
\STATE $cnt\leftarrow 0$
\REPEAT
\STATE $\mathcal F\leftarrow \mathcal D(M, \bar x, k)$
\STATE solve $\mathcal A(M, \bar x, \mathcal F)$ exactly and let $x^*$ denote the corresponding solution
\IF{$c^\top x^* < c^\top \bar x$}
\STATE $\bar x \leftarrow x^*$
\ELSE
\STATE $k\leftarrow \eta\cdot k$
\STATE $cnt\leftarrow cnt+1$
\ENDIF
\UNTIL {$cnt=C$ or time limit is reached}
\STATE \textbf{return} $\bar x$
\end{algorithmic}
\end{algorithm}
\end{minipage}\hfill
\begin{minipage}{0.510\textwidth}
\begin{algorithm}[H]
\caption{Two-Layer Large Neighborhood Search (TLNS)}
\begin{algorithmic}[1]
\label{alg:TLNS}
\STATE {\bfseries Input:} an MILP $M$,
   initial solution $\bar{x}$, fixing heuristic $\mathcal D$, presolve operator $\mathcal Q$, count limit $C$, neighborhood sizes $k_1, k_2$, and adaptive rates $\eta_1, \eta_2$.
\REPEAT
\STATE $\mathcal F \leftarrow \mathcal D(M, \bar x, k_1)$
\STATE $\{P, \bar y\}\leftarrow \mathcal Q(\mathcal A(M, \bar x ,\mathcal F), \bar x)$\hspace{0.1cm}\#presolve
\STATE $y^* \leftarrow \text{LNS}(P, \bar y, \mathcal D, C, k_2, \eta_2)$ \hspace{-0.3cm} 
\STATE $\{M, x^*\}\leftarrow \mathcal Q^{-1}(P, y^*)$ \hspace{0.43cm}\#postsolve
\IF {$c^\top x^* < c^\top \bar x$}
\STATE $\bar x\leftarrow x^*$
\ELSE 
\STATE $k_1\leftarrow \eta_1 \cdot k_1$
\ENDIF
\UNTIL {time limit is reached}
\STATE \textbf{return} $\bar x$
\end{algorithmic}
\end{algorithm}
\end{minipage}

The classic LNS method starts with an initial feasible solution and then gradually refines it by iteratively constructing a fixing neighborhood and optimizing the corresponding auxiliary problem.
The motivation for TLNS stems from the observation that auxiliary problems are still MILPs and hence can be further handled via LNS, rather than an off-the-shelf solver.
We outline the classic LNS in Algorithm~\ref{alg:LNS}.

\subsubsection{The Two-Layer Algorithm}

The TLNS framework can then be partitioned into two layers: an \textit{outer layer} and an \textit{inner layer}. 
Given an MILP (denoted by $M$) and an initial solution $\bar{x}$, we apply LNS iteratively to improve the incumbent solution, and this is called \enquote{outer layer}.
At each iteration, an auxiliary MILP has to be solved.
Again, we apply LNS to it, which we call \enquote{inner layer}.
The pseudocode of TLNS is provided in Algorithm~\ref{alg:TLNS}. 

\textbf{Outer Layer.}\quad 
Given an MILP $M$ and its initial solution $\bar{x}$, one can apply a fixing heuristic to build a neighborhood of size $k_1$, defining an auxiliary problem $\mathcal A(M, \bar x, \mathcal F)$.
Before we directly call LNS to solve this auxiliary problem, we need to deploy a critical operation: \textit{presolve}.
As~\cite{achterberg2020presolve} pointed out, presolve can be viewed as a collection of preprocessing techniques that reduce the size of and, more importantly, improve the \enquote{strength} of the given MILP, that is, the degree to which the constraints of
the formulation accurately describe the underlying polyhedron of integer-feasible solutions.
TLNS relies on a presolve operator $\mathcal Q$ to transform $M$ and the incumbent $\bar{x}$ into a reduced problem $P$ and $\bar{y}$, respectively.
The incumbent $\bar{x}$ will be improved from the inner layer and the neighborhood size $k_1$ increases if no improvements are made.
The outer layer will terminate if the time limit is reached. 

\textbf{Inner Layer.}\quad
The inner layer receives the presolved problem $P$ along with its feasible solution $\bar y$ from the outer layer.
A classic LNS is then invoked to optimize $P$.
Specifically, we employ a count limit as the stopping criterion for the inner layer LNS, as described in Algorithm~\ref{alg:LNS}. 
During each iteration, if the inner layer fails to find a better solution to $P$ within the neighborhood \text{i.e.} it is stuck in local minima, the neighborhood size $k_2$ should be increased to facilitate exploration of a broader search space.
Simultaneously, $cnt$ is incremented to keep track of the number of times the neighborhood size has been augmented.
In the end, a high-quality solution $y^*$ to $P$ is fed back to the outer layer where $\mathcal Q^{-1}$ transforms $y^*$ back to its counterpart $x^*$ in $M$.

\subsubsection{Comparison between LNS and TLNS}

Let LNS($\mathcal F$) denote a single LNS process with the fixing neighborhood defined by $\mathcal{F}$ (i,e. optimizing $\mathcal A(M, \bar x, \mathcal F)$).
For convenience, let LNS($\mathcal F^1:\mathcal F^H$) denote a process of applying LNS($\mathcal F^1$), LNS($\mathcal F^2$), ..., LNS($\mathcal F^H$) consecutively, with the superscript denoting its sequence. 
We utilize the subscript \enquote{1} and \enquote{2} to denote the outer and inner layer, respectively.
Let TLNS($\mathcal F_1, \mathcal F_2^1:\mathcal F_2^H$) denote a process of first enforcing the fixing neighborhood $\mathcal B(\mathcal F_1,\bar x)$ in the outer layer and then applying LNS($\mathcal F_2^1:\mathcal F_2^H$) in the inner layer.

\begin{remark}
\label{rmk:3.1}
    Compared to LNS($\mathcal F_1$), TLNS($\mathcal F_1,\mathcal F_2^1:\mathcal F_2^H$) exits search around the fixing neighborhood defined by $\mathcal F_1$ faster.
\end{remark}
Both methods result in the same auxiliary problem $\mathcal A(M, \bar x, \mathcal F_1)$.
The difference is that LNS($\mathcal F_1$) optimizes this problem via exact solvers while TLNS($\mathcal F_1,\mathcal F_2^1:\mathcal F_2^H$) addresses such an MILP via LNS($\mathcal F_2^1:\mathcal F_2^H$), a fast heuristic.
Firstly, exact solvers are based on the branch-and-bound framework enhanced with a dozen of modules that are critical to exactness but computationally expensive, such as cutting-planes, domain propagation and symmetry-breaking~\cite{wolsey1999integer}.
Alternatively, LNS($\mathcal F_2^1:\mathcal F_2^H$) channels all attention to searching for high-quality feasible solutions and is thus more efficient.
Secondly, TLNS($\mathcal F_1, \mathcal F_2^1:\mathcal F_2^H$) quickly identifies near-optimal solutions to $\mathcal A(M,\bar x, \mathcal F_1)$ and then exits search around the fixing neighborhood $\mathcal{B}(\mathcal{F}_1, \bar{x})$, moving towards a new neighborhood with potentially better solutions.

\begin{remark}
\label{rmk:3.2}
    Compared to LNS($\mathcal F_1\cup\mathcal F_2^1:\mathcal F_1\cup\mathcal F_2^H$), TLNS($\mathcal F_1,\mathcal F_2^1:\mathcal F_2^H$) saves presolving time.
\end{remark}

During step $h$, both methods utilize general-purpose solvers to optimize an auxiliary problem associated with $\mathcal F_1 \cup \mathcal F_2^h$.
For LNS($\mathcal F_1\cup\mathcal F_2^1:\mathcal F_1\cup\mathcal F_2^H$), the auxiliary problem is given as follows:  
\begin{equation}
\setlength\abovedisplayskip{1.5pt}
\setlength\belowdisplayskip{0.2pt}
    \begin{aligned}
          &  \underset{x \in S \cap \mathcal{B}\left(\mathcal F_1\cup\mathcal F_2^h, \bar{x}\right)}{\min} && c^\top x.
    \end{aligned} \label{eq:rmk2}
\end{equation}
Note that model~(\ref{eq:rmk2}) is the same as model~(\ref{eq:MILP}) except that some variable bounds are fixed.
As a result, presolving such models would become computationally costly, sometimes even exceeding its subsequent branch-and-bound tree search.
Let $T_p^1$ and $T_o$ denote the presolve time and the branch-and-bound search time, respectively.
Then the total used time for LNS($\mathcal F_1\cup\mathcal F_2^1:\mathcal F_1\cup\mathcal F_2^H$) is $H\times(T_p^1+T_o)$.
In TLNS($\mathcal F_1,\mathcal F_2^1:\mathcal F_2^H$), the auxiliary problem $\mathcal A(M, \bar x, \mathcal F_1)$ in the outer layer is presolved only once, producing an MILP of a reduced size (denoted by $P$).
Then an exact solver will be employed to first presolve $\mathcal A(P, \bar y, \mathcal F_2^h)$
(time $T_p^2$) with the subsequent branch-and-bound search (time $T_o$).
The total used time for TLNS($\mathcal F_1,\mathcal F_2^1:\mathcal F_2^H$) is $T_p^1+H\times(T_p^2+T_o)$.
Given that $T_p^2$ is much smaller than $T_p^1$, the saved time is $(H-1)\times T_p^1-H\times T_p^2$ when adopting TLNS($\mathcal F_1,\mathcal F_2^1:\mathcal F_2^H$).

The above two remarks elucidate the advantage of TLNS over LNS.
Typically, Remark~\ref{rmk:3.2} distinguishes between TLNS and LNS with smaller neighborhoods, highlighting the superiority of adopting such a nested approach over solely utilizing single-layer LNS with a small neighborhood.

\subsection{Learning-Enhanced TLNS}
We now utilize ML as a fixing heuristic in Algorithm~\ref{alg:TLNS}.
Let $s^t$ denote the state of an MILP $M$ with the incumbent solution $x^t$ in step $t$.
Our goal is to learn a policy $\pi_\theta(\cdot)$ that takes $s^t$ as the input and returns scores to determine the fixing neighborhood.
In the following, we first describe our training task and introduce the policy network, then we explain how we apply the learned policy for inference.

\subsubsection{Training}
Following previous works~\cite{huang2023searching,sonnerat2021learningnew}, we employ LB as the expert and collect samples for training.
\begin{definition}[LB neighborhood]
\label{def:LB_neighborhood}
Consider an MILP with $n$ variables, let $\bar{x}$ denote a feasible solution and  $k$ denote a distance cutoff parameter.
Then an \textit{LB neighborhood} is restricted to a ball around $\bar{x}$:
\begin{equation*}
    \begin{aligned}
        \mathcal{B}(k, \bar{x})   \coloneqq  \left\{x \in \mathbb{R}^n:  \lVert x - \bar{x} \rVert_1 \leq k  \right\},
    \end{aligned}
\end{equation*}
where $\lVert \cdot \rVert_1$ denotes $\ell_1$-norm.
\end{definition}
Using an LB neighborhood, we define a sub-MILP~(\ref{eq:LB_model}). 
\begin{equation}
\setlength\abovedisplayskip{1.5pt}
\setlength\belowdisplayskip{0.2pt}
    \begin{aligned}
        &  \underset{x \in S \cap \mathcal{B}\left(k, \bar{x}\right)}{\min} && c^\top x.
    \end{aligned}  \label{eq:LB_model}
\end{equation}

Model~(\ref{eq:LB_model}) is optimized by an MILP solver and let $x^*$ denote its optimal solution.
Since $x^* \in S$ and $c^\top x^* \leq c^\top \bar{x}$, $x^*$ becomes a new incumbent solution. 
Let a binary vector $a^*$ denote an action that can transform the previous incumbent $\bar{x}$ to the new one $x^*$, i.e., $a^*_i \coloneqq |\bar{x}_i - x^*_i|$.
We repeat the procedure of defining an LB neighborhood around the incumbent and optimizing the corresponding sub-MILP, until no objective improvement is achieved.

While optimizing a sub-MILP~(\ref{eq:LB_model}) in each iteration, we retrieve intermediate solution $\tilde{x}$ from the solution pool of an solver if $c^\top \bar x-c^\top \tilde x \geq \kappa_p (c^\top \bar x -c^\top x^*)$, with $0 < \kappa_p < 1$.
These solutions are not necessarily optimal but of high quality, defining a \textit{positive sample} set.
Specifically, let $\mathcal{S}_p$ denote such a set consisting of action vectors that can transform $\bar{x}$ to $\tilde{x}$.
We randomly flip elements of $a^*$ by $10\%$ while keeping the number of non-zero elements in $a^*$ unchanged, which generates a new action $a'$.
We then apply $a'$ to $\bar{x}$ and induce an MILP of form~(\ref{eq:model_auxi}) with $x'$ being the optimal solution. 
We accept action $a'$ as a \textit{negative sample} if $c^\top \bar x-c^\top  x' \leq \kappa_n (c^\top \bar x -c^\top x^*)$ with $0<\kappa_n\leq \kappa_p$.
Let $\mathcal{S}_n$ denote the set of negative samples.
Finally, let $\mathcal D \coloneqq \left\{(s,\mathcal S_p,\mathcal S_n)\right\}$ denote the set of collected data.

We utilize \textit{contrastive loss} for training. 
Formally, the loss function can be formulated as follows:
\begin{equation*}
    \begin{aligned}
            L(\theta) \coloneqq \sum\limits_{(s,\mathcal S_p, \mathcal S_n)\in \mathcal D} \frac{-1}{|\mathcal S_p|}\sum\limits_{a\in \mathcal S_p} \log \frac{\exp (a^\top \pi_\theta(s)/\tau)}{\sum\limits_{a'\in\mathcal S_n\cup\{a\}}\exp(a'^\top\pi_\theta(s)/\tau)},
    \end{aligned}
\end{equation*}
where $\tau$ is the temperature hyper-parameter.
The contrastive loss is deployed to bring $\pi_\theta(\cdot)$ closer to positive samples while simultaneously pushing it away from negative samples.
When $|\mathcal S_p|=1$ and $|\mathcal S_n|=0$ (e.g., only optimal solutions are kept as samples), we reduce the contrastive loss to the classic cross entropy loss used in~\cite{sonnerat2021learningnew}.

\subsubsection{Policy Network}
\label{sec:policy network}
The input of the policy is $s^t$ and the output $\pi_\theta(s^t)\in[0,1]^n$ assigns scores for each variable.
To encode $s^t$ based on bipartite graph representations, previous works~\cite{huang2023searching,sonnerat2021learningnew} adopted a rich set of features including information derived from solving linear programs (LPs).
However, given the computational demands of these LP-based features in large-scale problems, we chose to include only those features that can be efficiently computed, as is deployed in~\cite{pns}.

Regarding the network architecture, while GNNs are naturally suited for bipartite graphs, the limited depth restricts interactions between distant variable nodes.
Besides, the classic self-attention~\cite{vaswani2017attention} layer with $O(n^2)$ complexity could be computationally prohibitive for large-scale MILPs.
To address these limitations, we employ the Simplified Graph Transformer~\cite{wu2024simplifying}, which expands the GNNs' receptive field through a global attention module while maintaining a lightweight structure due to its linear attention mechanism.
Specifically, we first employ the attention module described in model~(\ref{eq:sgt}) to aggregate information across the entire graph.
The output then serves as the initial node embedding of the subsequent GNN module, where we incorporate two interleaved half-convolution layers~\cite{gasse2019exact}.
Finally, the embeddings of variable nodes are transformed into scalars within $[0,1]$ through 2-layer perceptrons alongside a sigmoid function.

\subsubsection{Inference} 
In LNS, when building a fixing neighborhood of size $r$, we apply the learned policy $\pi_\theta(\cdot)$ to inform the neighborhood.
We employ a sampling strategy to randomly select $r$ variables to be unfixed without replacement according to $\pi_\theta(s)$, where variables with higher scores are more likely to be selected.

\section{Numerical Experiments} 
\label{sec:experiments}
In this section, we evaluate the performance of our proposed algorithm and compare it with other methods. The code will be made publicly available upon publication.

    \begin{table}[t]
\caption{Average size of each benchmark instance, the SMALL instances are used for data collection and training and the LARGE instances are used for testing}
\label{table:instances}
\centering
\resizebox{\columnwidth}{!}{\begin{tabular}{ccccccccc}
\toprule
\multirow{2}{*}{}& \multicolumn{4}{c}{SMALL} & \multicolumn{4}{c}{LARGE}\\
\cmidrule(lr){2-5}\cmidrule(lr){6-9}
&SC & CA & MIS & MVC & SC & CA & MIS & MVC\\

\midrule
\text{\# variables}& 4,000& 4,000&6,000 & 1,000& 16,000& {100,000} &100,000 &20,000\\
\text{\# constraints}& 5,000& 2,662&15,157 & 65,100&20,000 & {794,555} &5,001,669 &3,960,000\\
\text{\# non-zeros}& 1,000,000& 22,757& 30,314& 130,200&16,000,000 & 4,000,000 & 10,003,338 &7,920,000\\
\bottomrule
\end{tabular}}
\end{table}
\subsection{Setup}
\label{sec:4.1_experiment_setup}

\subsubsection{Benchmarks}
In our evaluation, we assess our algorithm on four widely used $\mathcal{NP}$-hard problem benchmarks—Set Cover (SC), Combinatorial Auction (CA), Maximum Independent Set (MIS), and Minimum Vertex Cover (MVC)—following the instance generation procedures of~\cite{huang2023searching}. 
Each SC instance is generated with $5,000$ items and $4,000$ subsets, while every CA instance is generated with $4,000$ bids and $2,000$ items. 
MIS and MVC are graph-related problems and they are generated from random graphs with $6,000$ and $1,000$ nodes, average degrees of $5$ and $130$, respectively.
For each benchmark, we create $1,000$ instances (denoted by SMALL), split into training and validation sets of $900$ and $100$ instances, respectively. 
Additionally, we generate $20$ larger instances (denoted by LARGE) for each benchmark, split in half for test and validation sets.
LARGE instances range from 4 to 25 times the size of SMALL instances with MILP sizes specified in Table~\ref{table:instances}. 
These MILPs have up to $100$ thousand variables and $5$ million constraints, becoming computationally prohibitive for general-purpose solvers.

\subsubsection{Evaluation Configurations}
All evaluations are performed under the same configuration.
The evaluation machines include 12th Gen Intel(R) Core(TM) i9-12900K CPUs with Nvidia GeForce RTX 3090 GPUs. 
For off-the-shelf MILP solvers, Gurobi $10.0.2$~\cite{gurobi} and SCIP $8.0.4$~\cite{scip}
are utilized in our experiments. 
The time limit for running each experiment
is set to $1,000$ seconds since a tail-off of solution qualities was often observed after that.

\subsubsection{Data collection \& Training}
For each training instance, we utilize Gurobi with a solution limit of $1$ to generate the very first incumbent solution.
Using this incumbent, we apply the LB heuristic with a fine-tuned neighborhood size of $100$, $400$, $500$, and $75$ for SC, CA, MIS, and MVC, respectively.
The resulting LB-MILP is optimized by Gurobi with a time limit of $1,500$ seconds.
We employ the contrastive loss with $\tau$ being equal to $0.07$.
The batch size is set to $32$ and Adam~\cite{kingma2014adam} with a learning rate of $0.001$ is utilized as the optimizer.
We remark that in our experiments: (i) The Graph Transformer models are trained using SMALL instances but applied to LARGE ones; (ii) the same models are deployed in LNS as well as both layers of TLNS.

\subsubsection{Metrics}
In order to assess the performance of different methods, we employ two metrics: (i) \textit{primal bound} (PB), which refers to the objective value; (ii) \textit{primal integral} (PI), which measures the integral of primal gap with respect to runtime, where the primal gap denotes a normalized difference between the primal bound and a pre-specified best known objective value.
The best known objectives correspond to the best solutions returned by all methods evaluated in our experiments with a longer time limit of $3,600$ seconds.
The PB value demonstrates the quality of returned solutions at runtime while the PI value provides insight into the speed of identifying better solutions.
Note that we choose not to report the primal gap as a measure of optimality in our experiments since for large-sized instances, both Gurobi and SCIP could neither identify optimal/near-optimal solutions nor provide tight dual bounds with a reasonable time limit (e.g. $500,000$ seconds).

Since all four benchmarks considered in our experiments entail minimization problems, smaller PB and PI values imply better computational performances.

\begin{wraptable}{r}{6cm}
\vskip -0.3in
\caption{Neighborhood sizes}
\label{table:neighborhood size}
\centering
\resizebox{0.5\columnwidth}{!}{
\begin{tabular}{ccccc}
\toprule
{Dataset}& \texttt{R-LNS} & \texttt{R-TLNS} & \texttt{CL-LNS} & \texttt{CL-TLNS} \\
\midrule

{SC} & 4,000 & 8,000/1600  & 175 & 500/120  \\

CA & {35,000} &{60,000/3,000} & {35,000} & {60,000/3,000}\\

MIS & 40,000 & 70,000/7,000 & 12,500 & 30,000/7,000\\

MVC & 10,000& 15,000/1,250& 1,250 & 5,500/1,000 \\

\bottomrule
\end{tabular}
}
\vskip -0.2in
\end{wraptable}  
To showcase the effectiveness of TLNS, we conduct the following progressive experiments: (i) comparing TLNS with LNS under classic settings~(Section~\ref{sec:comparison_tradition}); (ii) comparing TLNS with LNS under learning-based settings~(Section~\ref{sec:comparison_learning}); (iii) comparing against state-of-the-art MILP solvers~(Section~\ref{sec:comparison_solver}); and (iv) comparing TLNS against LP-free heuristics (Section~\ref{sec:vs heur}).
The neighborhood sizes for each method are fine-tuned individually and summarized in Table~\ref{table:neighborhood size}. For TLNS, the neighborhood sizes (e.g., $500/120$) signify that $500$ variables remain unfixed in the outer layer, while $120$ variables remain unfixed in the inner layer.
Readers are referred to Section~\ref{sec:ablation} for ablation experiments.

\subsection{Comparison between TLNS and LNS (classic)} 
\label{sec:comparison_tradition}

\begin{figure}
\centering
\includegraphics[width=0.5\linewidth]{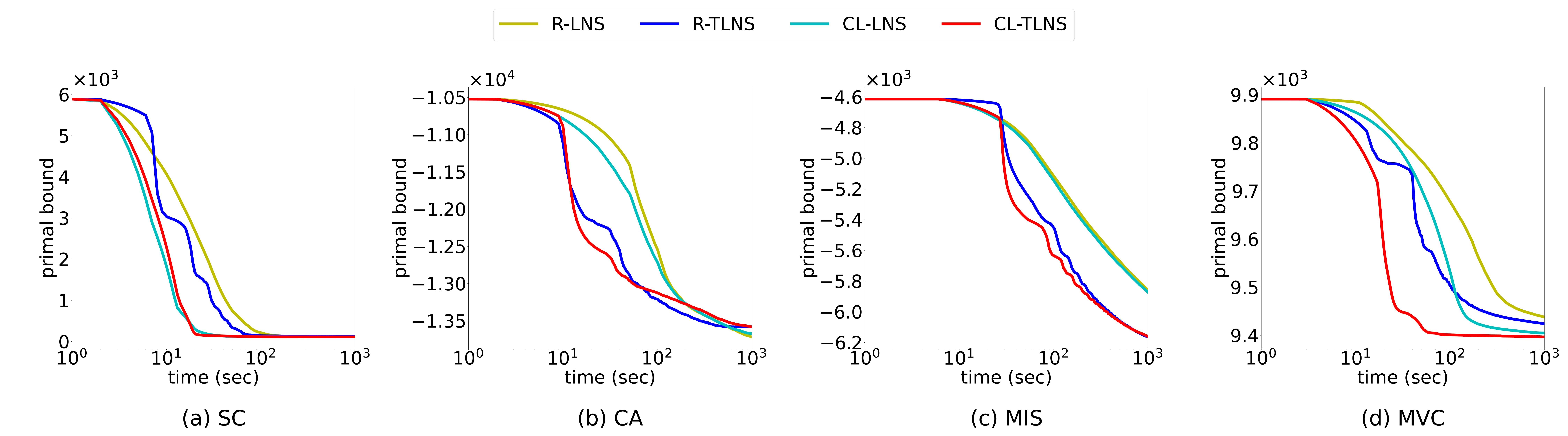}

\centering
\subfigure[SC]{
\label{SC}
\begin{minipage}[b]{0.23\linewidth}
\includegraphics[width=1\linewidth]{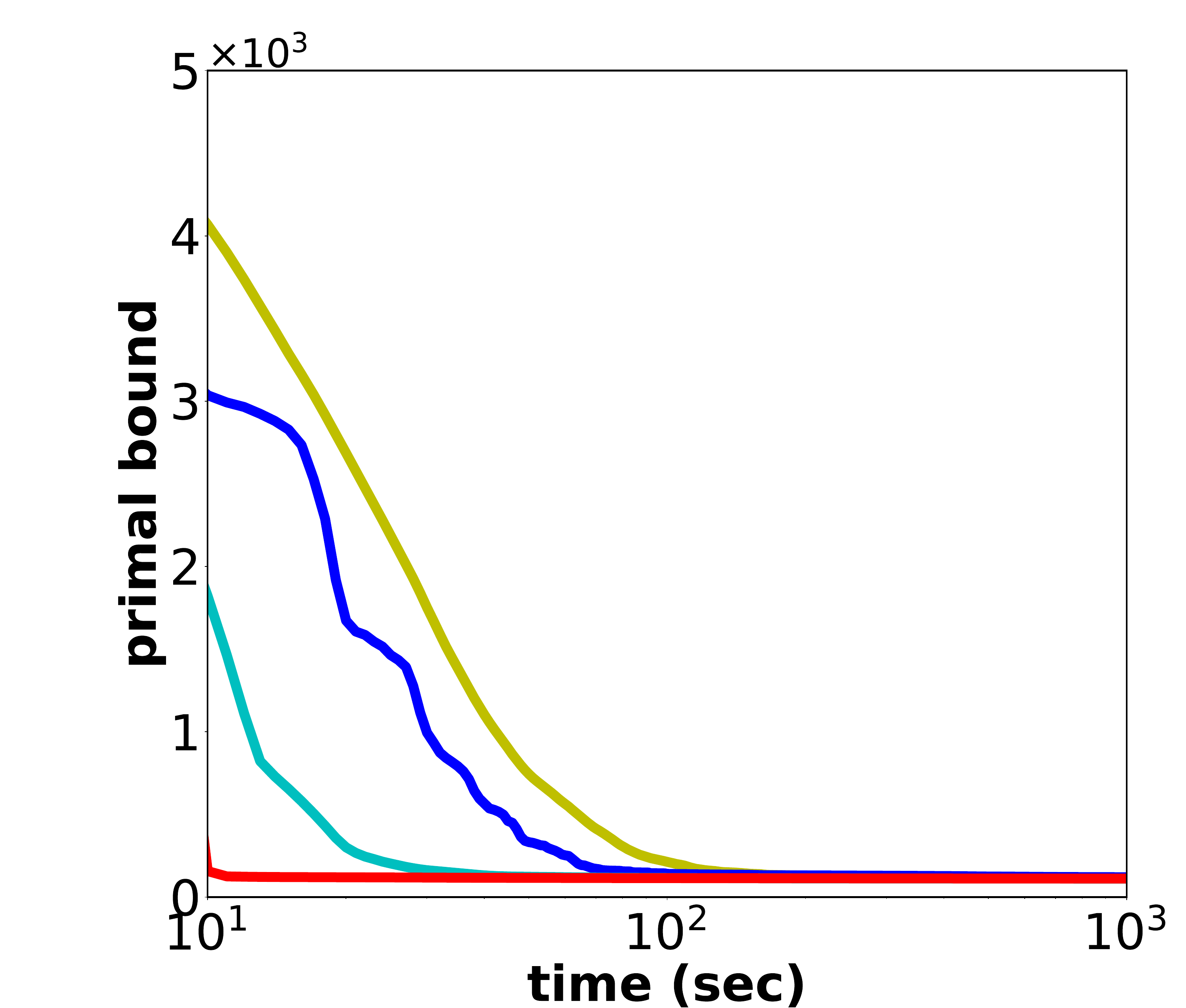}
\end{minipage}}
\subfigure[CA]{
\label{CA}
\begin{minipage}[b]{0.23\linewidth}
\includegraphics[width=1\linewidth]{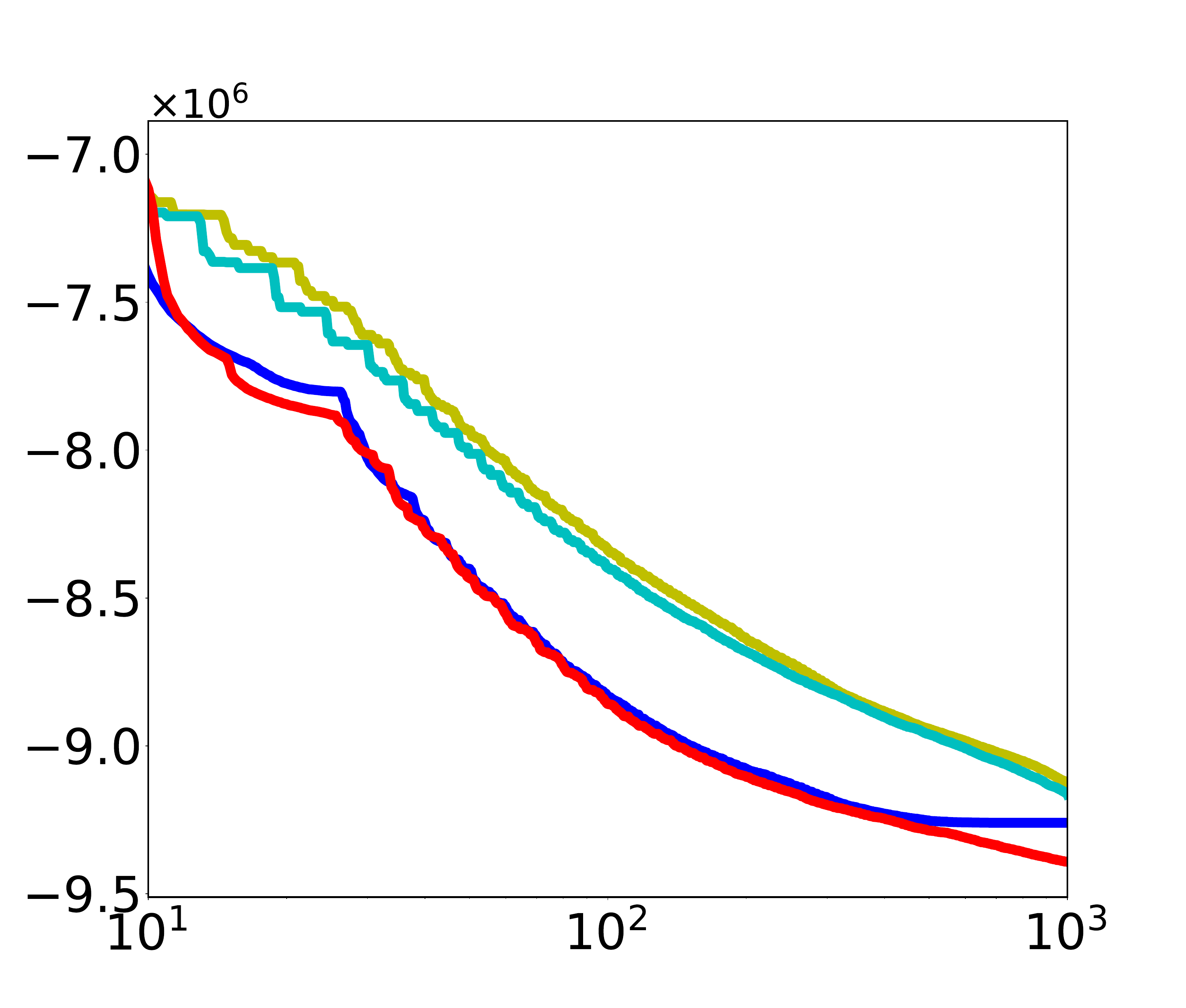}
\end{minipage}}
\subfigure[MIS]{
\label{MIS}
\begin{minipage}[b]{0.23\linewidth}
\includegraphics[width=1\linewidth]{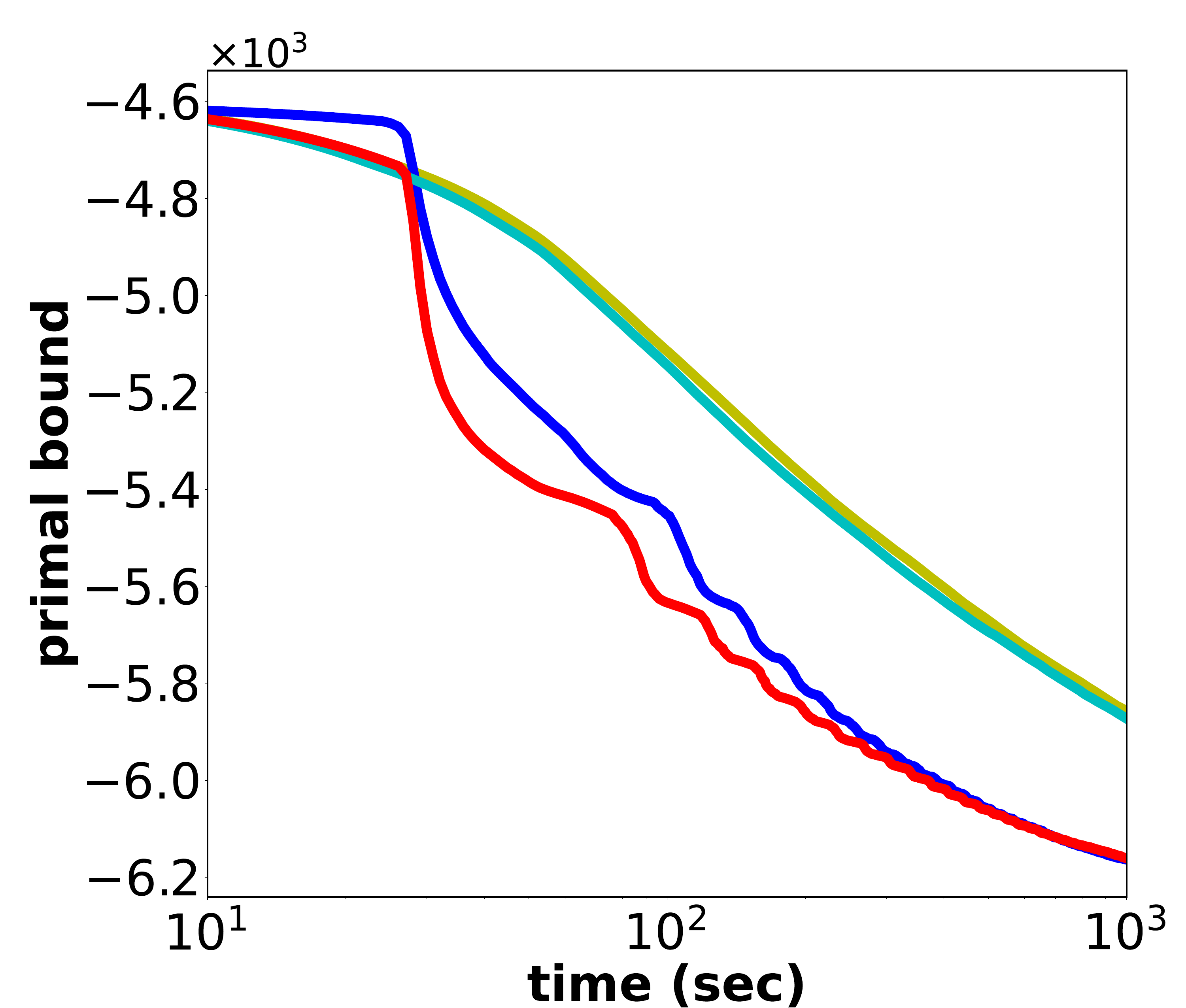}
\end{minipage}}
\subfigure[MVC]{
\label{MVC}
\begin{minipage}[b]{0.23\linewidth}

\includegraphics[width=1\linewidth]{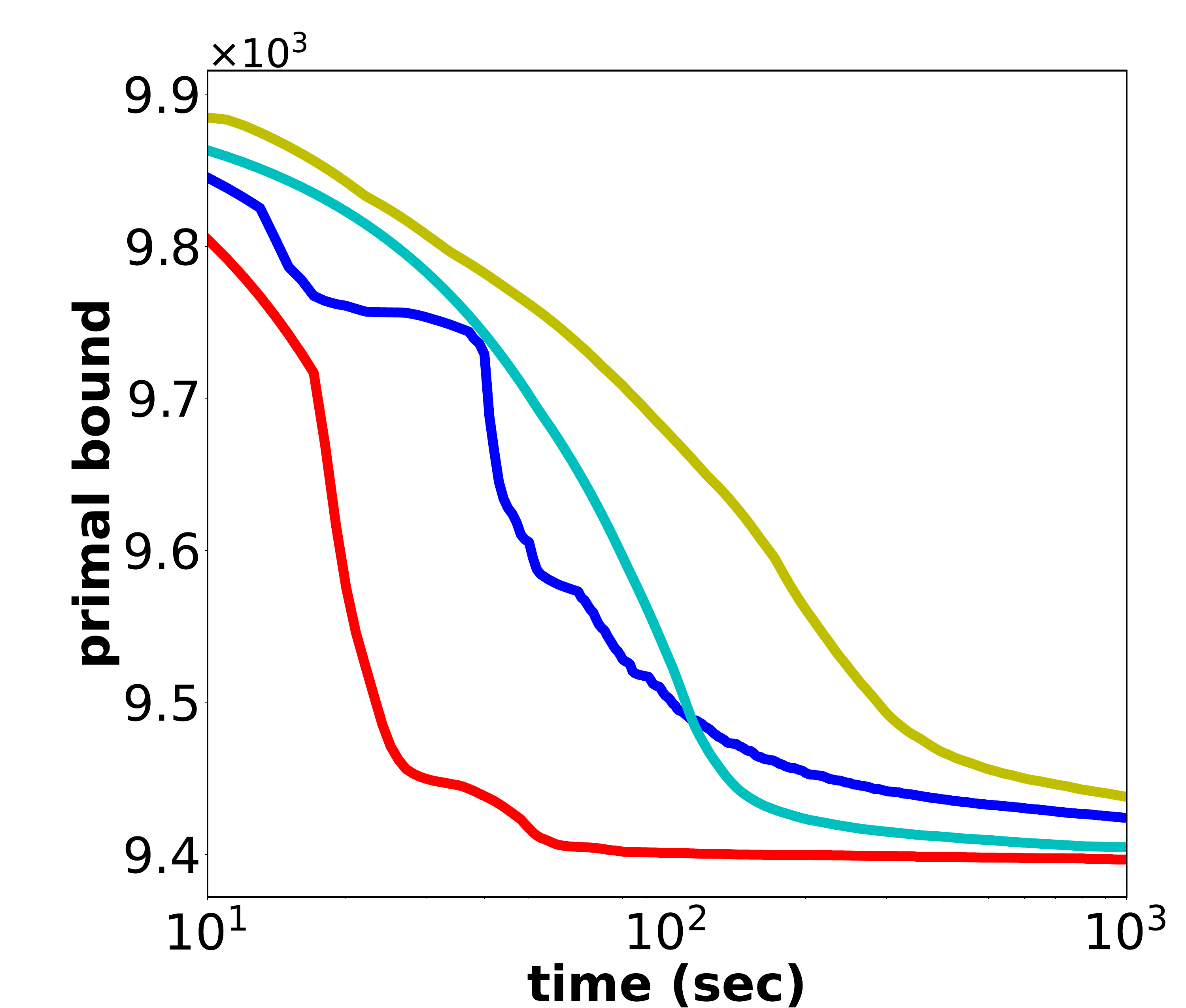}
\end{minipage}}
\caption{The PB, as a function of runtime, averaged over 10 instances. 
Lower PB values imply better performance.}
\label{figure:pb-runtime}
\vskip -0.2 in
\end{figure}

We hypothesize that our proposed TLNS outperforms LNS irrespective of neighborhood-fixing heuristics.
To demonstrate this, we first compare TLNS against LNS under classic learning-free settings.
Among various heuristics for fixing variables, we employ a \textit{random} heuristic, which has been used in~\cite{rothberg2007evolutionary} and is straightforward to implement as it selects variables randomly.
We denote LNS and TLNS with the random heuristic by \texttt{R-LNS} and \texttt{R-TLNS}, respectively.
The very first incumbent solution is again provided by Gurobi with a solution limit of one.
We employ SCIP as the off-the-shelf MILP solver in both \texttt{R-LNS} and \texttt{R-TLNS} since the presolve operator in Gurobi is inaccessible.

Fig~\ref{figure:pb-runtime} depicts the average PB as a function of runtime.
Clearly, the PB value of \texttt{R-TLNS} (blue) is smaller than that of \texttt{R-LNS} (yellow) for most of the runtime, demonstrating that our two-layer neighborhood search method indeed exceeds the single-layer version. 
Table~\ref{table:TLNS_vs_LNS} presents PI and PB values at a time limit of $1,000$ seconds, averaged over $10$ instances for each dataset, along with their respective standard deviations.
Bold values in columns \enquote{\texttt{R-LNS}} and \enquote{\texttt{R-TLNS}} indicate superior performances.
Column \enquote{Gain} presents the improvement of \texttt{R-TLNS} over \texttt{R-LNS} for each metric.
For instance, in the case of MVC, the PI values for \texttt{R-LNS} and \texttt{R-TLNS} are $14.5$ and $9.0$ respectively.
The gain is computed as $(14.5-9.0)/14.5 = 37.9\%$.
Notably, in terms of PI, \texttt{R-TLNS} surpasses \texttt{R-LNS} across all datasets, particularly exhibiting a remarkable improvement of {$51.3\%$ on the CA dataset.}
This indicates that \texttt{R-TLNS} produces high-quality solutions faster than \texttt{R-LNS}.
In terms of PB at $1,000$ seconds, \texttt{R-TLNS} produces better solutions to {CA and MIS}, equivalently good ones to {MVC}, and slightly worse ones to SC, compared with \texttt{R-LNS}.
The latter phenomenon is potentially due to the fact that both methods stagnate in the later runtime but \texttt{R-LNS} relying on exact solvers is capable of producing optimal solutions to auxiliary problems.
We claim that \texttt{R-TLNS} generally outperforms \texttt{R-LNS}.










\begin{table}[b]
    \caption{PI and PB values at $1,000$ seconds for \texttt{R-LNS}, \texttt{R-TLNS}, \texttt{CL-LNS} and \texttt{CL-TLNS} averaged over $10$ instances for each benchmark, along with their standard deviations. 
    Lower PI/PB values imply better performances.}
    \label{table:TLNS_vs_LNS}
    \renewcommand\arraystretch{1.3}
    \resizebox{\columnwidth}{!}{
    \begin{tabular}{cccccccc}
    \toprule
    \multicolumn{2}{c}{Dataset}& {\texttt{R-LNS}} & {\texttt{R-TLNS}} & Gain & {\texttt{CL-LNS}} & {\texttt{CL-TLNS}} & Gain \\
    \cmidrule(lr){1-2}\cmidrule(lr){3-5}\cmidrule(lr){6-8}
    \multirow{2}{*}{SC} & PI & 1695.7$\pm$178.7 & \textbf{1250.9$\pm$137.8} &\textbf{26.2\%}  & 871.8$\pm$58.5 & \textbf{633.5$\pm$ 43.6} & \textbf{27.3\%} \\
                        & PB & \textbf{118.0$\pm$3.8} & 121.3$\pm$6.17 & -2.7\% & \textbf{112.7$\pm$2.7} & 113.0$\pm$2.9&-0.2\% \\
    \cmidrule(lr){1-2}\cmidrule(lr){3-5}\cmidrule(lr){6-8}
    
    \multirow{2}{*}{CA} & PI & 64.2$\pm$3.4 & \textbf{31.2$\pm$4.3} & \textbf{51.3\%} & 60.5$\pm$2.3 & \textbf{25.5$\pm$1.7} &\textbf{57.8\%}  \\
                        & PB & {-9124241.2$\pm$29637.7} & \textbf{-9260001.4$\pm$34248.5} & \textbf{1.4\%} & -9162882.8$\pm$38388.0 & \textbf{-9391997.4$\pm$50556.4} &\textbf{2.5\%}  \\
    \cmidrule(lr){1-2}\cmidrule(lr){3-5}\cmidrule(lr){6-8}
    
    \multirow{2}{*}{MIS} & PI & 158.5$\pm$4.5 & \textbf{92.8$\pm$2.1} & \textbf{41.4\%}  &  149.4$\pm$17.0 & \textbf{50.6$\pm$ 4.9} & \textbf{66.1\%} \\
                        & PB & -5740.5$\pm$22.7 & \textbf{-6134.7$\pm$17.6} & \textbf{6.8\%} & -5897.8$\pm$85.6 & \textbf{-6435.8$\pm$16.8} &\textbf{9.1\%}  \\
    \cmidrule(lr){1-2}\cmidrule(lr){3-5}\cmidrule(lr){6-8}
    \multirow{2}{*}{MVC} & PI & 14.5$\pm$0.8 & \textbf{9.0$\pm$0.7} & \textbf{37.9\%} & 5.1$\pm$0.6 & \textbf{3.7$\pm$ 0.3} & \textbf{26.5\%} \\
                        & PB & 9446.6$\pm$39.9 & \textbf{9442.7$\pm$ 42.5}&\textbf{0.04\%} & 9401.8$\pm$39.6 & \textbf{9394.8$\pm$41.5} &\textbf{0.07\%} \\
    \bottomrule
    \end{tabular}}
\vskip -0.2in
\end{table}

\subsection{Comparison between TLNS and LNS (learning)}
\label{sec:comparison_learning}
We proceed to show that TLNS performs better than LNS in learning-based settings.
The very recent effort of enhancing LNS with ML is the work of~\cite{huang2023searching}, where the authors adopted contrastive learning to build fixing neighborhoods.
We apply this technique in the TLNS algorithm, denoting the two methods as \texttt{CL-LNS} and \texttt{CL-TLNS}, respectively. \texttt{SCIP} is used as the off-the-shelf MILP solver within both methods. \texttt{CL-LNS} and \texttt{CL-TLNS} are then evaluated on four benchmark datasets.
To ensure fairness, the same trained models are utilized in both methods.

We plot PB as a function of runtime in Fig~\ref{figure:pb-runtime}.
The PB curves of \texttt{CL-TLNS} (red) are almost consistently below that of \texttt{CL-LNS} (cyan).
Table~\ref{table:TLNS_vs_LNS} presents the PI/PB metrics of both methods at the time limit of $1,000$ seconds.
In terms of PI, \texttt{CL-TLNS} significantly surpasses \texttt{CL-LNS} across all four benchmark datasets--with an improvement ranging from $26.5\%$ to $66.1\%$.
In terms of PB, \texttt{CL-TLNS} produces better solutions to CA and MIS, equivalently good ones to MVC, and slightly worse ones to SC, compared with \texttt{CL-LNS}.
We found out that MVC instances were solved to near-optimality, hence \texttt{CL-TLNS} and \texttt{CL-LNS} achieve comparable PB metrics. 
As for SC, the mildly worse performance of \texttt{CL-TLNS} can be again attributed to stagnation, as discussed in Section~\ref{sec:comparison_tradition}.
We can claim that \texttt{CL-TLNS} generally outperforms \texttt{CL-LNS}.

We now compare the learning-guided fixing method with a random heuristic used in TLNS.
From Fig~\ref{figure:pb-runtime}, the PB value of \texttt{CL-TLNS} (red) is almost consistently lower than that of \texttt{R-TLNS} (blue) on all four benchmark datasets.
Notably, switching from a random heuristic to a learning-based one benefits TLNS significantly {across all datasets.}
From Table~\ref{table:TLNS_vs_LNS}, we compare both PB and PI metrics of \texttt{CL-TLNS} with those of \texttt{R-TLNS}. 
Again, we can claim that \texttt{CL-TLNS} generally performs better than \texttt{R-TLNS}.

\subsection{Comparison against MILP Solvers}
\label{sec:comparison_solver}
\begin{table}[b]
\caption{{PI and PB at $1,000$ seconds for SCIP, Gurobi, \texttt{GRB(NoRel)} and \texttt{CL-LNS(NoRel)}, averaged over $10$ instances for each benchmark, along with their standard deviations. 
Lower PI/PB values imply better performances.
}}
\label{table:vs solver}
\renewcommand\arraystretch{1.3}
\resizebox{\columnwidth}{!}{
\begin{tabular}{ccccccccccc}

\toprule
\multicolumn{2}{c}{Dataset}& {CL-TLNS} & {SCIP} & {Gurobi} & Gain & {\texttt{GRB(NoRel)}} & {\texttt{CL-LNS(NoRel)}} & Gain \\
\cmidrule(lr){1-3}\cmidrule(lr){4-6}\cmidrule(lr){7-9}

\multirow{2}{*}{SC} 
  & PI& 633.5$\pm$43.6&262487$\pm$7482 & \underline{16294$\pm$1035}& {96.1\%}& 15378.4$\pm$1028.7 & \underline{2473.6$\pm$242.5} &74.4\%\\
& PB & 113.0$\pm$2.9&  \underline{116.6$\pm$2.6} & 120.9$\pm$2.7&{3.0\% }&  114.9$\pm$2.7 &\underline{112.6$\pm$2.8} &-0.3\% \\

\cmidrule(lr){1-3}\cmidrule(lr){4-6}\cmidrule(lr){7-9}
\multirow{2}{*}{CA} 
  & PI&25.5$\pm$1.7&285.3$\pm$5.7& \underline{180.3$\pm$3.4}& {85.8\%}&\underline{41.3$\pm$4.9}& {55.9$\pm$4.6} & 38.2\%\\
& PB&-9391997.4$\pm$50556.4 & -7071993.8$\pm$50583.6& \underline{-7728867.4$\pm$36376.9} &{21.5\%}& -9283749.8$\pm$57926.9& \underline{-9300736.2$\pm$33056.6} & 1.0\%\\

\cmidrule(lr){1-3}\cmidrule(lr){4-6}\cmidrule(lr){7-9}
\multirow{2}{*}{MIS} 
& PI& 50.6$\pm$4.9&--& \underline{284.9$\pm$3.9}&82.2\%&182.5$\pm$5.7 &\underline{68.0$\pm$7.6}&25.6\%\\
&PB&-6435.8$\pm$16.8 &--&\underline{-4606.7$\pm$20.2} &39.7\%&-6009.2$\pm$28.4& \underline{-6298.4$\pm$17.7}&2.1\%\\
 
\cmidrule(lr){1-3}\cmidrule(lr){4-6}\cmidrule(lr){7-9}
\multirow{2}{*}{MVC} 
& PI&3.7$\pm$0.3 & {45.3$\pm$0.8} & \underline{41.5$\pm$3.3} &90.9\%&16.2$\pm$0.3 & \underline{3.7$\pm$0.4} &0.0\%\\
&PB& 9394.8$\pm$41.5 & \underline{9602.8$\pm$42.4} &9751.1$\pm$50.2&2.1\%&9423.6$\pm$38.3 & \underline{9394.7$\pm$39.4}& 0.0\%\\

\bottomrule

\end{tabular}}

\end{table}  
We compare \texttt{CL-TLNS} against Gurobi and SCIP across four datasets. 
To ensure a fair comparison between exact solvers and heuristics, we enforce both solvers to apply their internal heuristics more aggressively. 
In particular, we use \enquote{model.setHeuristics(SCIP\_PARAMSETTING.AGGRESSIVE)} for SCIP and set the parameter \enquote{MIPFocus = 1} for Gurobi.
Table~\ref{table:vs solver} exhibits the averaged PI and PB metrics.
Note that \enquote{--} in Column \enquote{SCIP} indicates that SCIP is incapable of handling MIS problems due to memory limits.
The underlined values in columns \enquote{SCIP} and \enquote{Gurobi} signify better performances between these two solvers while column \enquote{Gain} represents the improvement of \texttt{CL-TLNS} over the superior solver in terms of respective metrics.
The computational results show that \texttt{CL-TLNS} consistently outperforms Gurobi and SCIP across all benchmarks, achieving an improvement of up to $96.1\%$ and $39.7\%$ in PI and PB, respectively.

\subsection{Comparison against LP-free Heuristics}
\label{sec:vs heur}

In Section~\ref{sec:comparison_solver}, we compare our approach with advanced MILP solvers. However, for large-scale problems, these solvers can be inefficient due to their exact search nature and the high computational cost of solving LP relaxations, which do not necessarily contribute to improving PB values. To address these inefficiencies, the current section introduces an alternative baseline by employing LP-free heuristics. Specifically, we use Gurobi with \enquote{NoRel} (No Relaxation) heuristics to (i) directly solve MILPs (denoted as \texttt{GRB(NoRel)}) and (ii) serve as the underlying solver for \texttt{CL-LNS} (denoted as \texttt{CL-LNS(NoRel)}).
We then compare \texttt{CL-TLNS} with both \texttt{GRB(NoRel)} and \texttt{CL-LNS(NoRel)}, with the results presented in Table~\ref{table:vs solver}. The underlined values in the columns labeled \enquote{\texttt{GRB(NoRel)}} and \enquote{\texttt{CL-LNS(NoRel)}} indicate the better performance between these two approaches, while the \enquote{Gain} column reflects the improvement of \texttt{CL-TLNS} over the superior approach in terms of the respective metrics.
In terms of PI, \texttt{CL-TLNS} significantly outperforms the other baselines on SC, CA, and MIS, while achieving comparable results to \texttt{CL-LNS(NoRel)} on MVC. Regarding PB, \texttt{CL-TLNS} delivers the best solutions for CA and MIS but slightly underperforms \texttt{CL-LNS(NoRel)} on SC, which can be attributed to stagnation. Notably, although \texttt{CL-LNS(NoRel)} relies on Gurobi as its underlying solver, \texttt{CL-TLNS} still outperforms it, despite using SCIP as its MILP solver. We conclude that \texttt{CL-TLNS} outperforms both \texttt{GRB(NoRel)} and \texttt{CL-LNS(NoRel)}.

It is also noteworthy that, among the three heuristics, the number of LNS layers in \texttt{GRB(NoRel)}, \texttt{CL-LNS(NoRel)}, and \texttt{CL-TLNS} are 0, 1, and 2, respectively. As the number of layers increases, performance improves. This indicates that adding more layers can be an effective strategy for tackling large-scale problems.

\subsection{Ablation Study}
\label{sec:ablation}
This section details ablation experiments comparing various model architectures.  
Specifically, SGT refers to the simplified graph transformer model introduced in Section~\ref{sec:policy network}. 
GCN denotes the classic graph neural network employing half-convolutions (equivalent to SGT without the attention layers), while GAT represents the Graph Attention Network with half-convolutions, as implemented in \cite{huang2023searching}. 
For all models, the number of layers and hidden dimensions are uniformly set to $2$ and $32$, respectively.

\begin{table}[htpb]
\caption{PI and PB at 1,000 for CL-TLNS with model architecture SGT, GCN, GAT. Lower PI/PB values imply better performances.}
\label{talbe:ablation}
\centering
\tabcolsep=1mm
\begin{tabular}{ccccc}
\toprule
\multicolumn{2}{c}{Dataset}& {SGT} & {GCN} & {GAT}\\
\midrule

\multirow{2}{*}{SC} 
& PI& \textbf{633.5$\pm$43.6} & 640.6$\pm$33.8 & {OOM}\\
& PB & \textbf{113.0$\pm$2.9} &  113.1$\pm$4.0  & {OOM}\\

\midrule

\multirow{2}{*}{CA} 
& PI& \textbf{25.5$\pm$1.7} & {196.7$\pm$6.6} & 35.1$\pm$3.5\\
& PB & \textbf{-9391997.4$\pm$ 50556.4}& {-8035948.2$\pm$82793.1} & {-9358284.7$\pm$16573.2}\\
  
 \midrule

\multirow{2}{*}{MIS} 
& PI& \textbf{50.6$\pm$4.9} & 170.3$\pm$3.3 & 78.2$\pm$3.5 \\
&PB & \textbf{-6435.8$\pm$16.8} & -5602.9$\pm$37.8  & -6162.3$\pm$13.7\\

\midrule
\multirow{2}{*}{MVC} 
& PI& 3.7$\pm$0.3& {5.25$\pm$0.5} & \textbf{3.17$\pm$0.3} \\
&PB & 9394.8$\pm$41.5 & {9401.1$\pm$42.2} & \textbf{9394.1$\pm$41.8}\\

\bottomrule
\end{tabular}
\end{table}  

The results are exhibited in Table~\ref{talbe:ablation}.
On the one hand, our SGT model surpasses GCN, leveraging its expanded receptive fields to enhance performance. On the other hand, SGT outperforms GAT on the CA and MIS tasks while delivering comparable results on MVC. Notably, GAT is significantly more memory-intensive and can encounter memory-related issues (e.g., out-of-memory (OOM) errors on SC), which SGT manages to avoid.

\section{Conclusions}
\label{sec:conclusion}
This paper proposes a learning-enhanced TLNS method especially for addressing large-scale MILPs.
Classic LNS methods refine incumbent solutions by building a particular neighborhood and searching within such a region by optimizing auxiliary MILPs via off-the-shelf solvers while our proposed TLNS goes one step further and solves auxiliary problems via LNS.
Graph transformer models are incorporated into TLNS for guiding neighborhood construction, boosting the performance of TLNS.
We argue that learning-based TLNS would outperform classic LNS and demonstrate this in our experiments.
The results show that TLNS achieved significantly better performances in identifying high-quality solutions within a short time frame.
An immediate research direction is to generalize TLNS and to extend it to the \textit{multi-layer} LNS, where LNS is applied recursively to solve subsequent sub-MILPs.
The \textit{multi-layer} LNS is promising to address extremely large-scale MILPs while reducing reliance on off-the-shelf solvers.

\bibliographystyle{splncs04}
\bibliography{ref}

\end{document}